\documentclass{gtart_h}  

\def\ifplaintex{\expandafter\ifx\csname documentclass\endcsname\relax}

\ifplaintex 
\else
\fi

\expandafter\ifx\csname beginpicture\endcsname\relax
\expandafter\ifx\csname documentclass\endcsname\relax
\input pictex \else
\input prepictex \input pictex \input postpictex \fi\fi

\def\gt{{\mathsurround=0pt\it $\cal G\mskip-2mu$eometry \&\ 
$\cal T\!\!$opology}}        

\def\gtp{{\mathsurround=0pt\it $\cal G\mskip-2mu$eometry \&\ 
$\cal T\!\!$opology $\cal P\!$ublications}}  

\def\lognumber#1{\def\thelognumber{#1}}
\def\volumenumber#1{\def\thevolumenumber{#1}}
\def\papernumber#1{\def\thepapernumber{#1}}
\def\volumeyear#1{\def\thevolumeyear{#1}}

\def\pagenumbers#1#2{\def\startpage{#1}\def\finishpage{#2}}
\def\published#1{\def\publishdate{#1}}
\def\proposed#1{\def\theproposer{#1}}
\def\seconded#1{\def\theseconders{#1}}
\def\received#1{\def\receiveddate{#1}}
\def\revised#1{\def\reviseddate{#1}}
\def\accepted#1{\def\accepteddate{#1}}
\def\asciititle#1{\def\theasciititle{#1}}

\long\def\asciiabstract#1{\long\def\theasciiabstract{#1}}
\def\asciikeywords#1{\def\theasciikeywords{#1}}

\let\\\par\let\thelognumber\relax
\let\thevolumenumber\relax\let\thepapernumber\relax
\let\thevolumeyear\relax\let\thesamplenumber\relax\let\startpage\relax
\let\finishpage\relax\let\publishdate\relax\let\receiveddate\relax
\let\reviseddate\relax\let\accepteddate\relax\let\theasciititle\relax
\let\theasciiauthors\relax
\let\theasciiabstract\relax\let\theasciikeywords\relax
\let\theasciiemail\relax\let\theshortauthors\relax\let\theshorttitle\relax

\long\def\maketitlep{   

\count0=\startpage

\gt\hfill      
\beginpicture
\setcoordinatesystem units <0.33truein, 0.33truein> point at 2.2 0.9
\setplotsymbol ({$\cal G$})
\plotsymbolspacing=9truept
\circulararc 315 degrees from 0 1 center at 0 0
\setplotsymbol ({$\cal T$})
\circulararc 315 degrees from 1 -1 center at 1 0
\endpicture
\break
{\small\ifx\thesamplenumber\relax 
Volume \else Sample
\fi\thevolumenumber\ (\thevolumeyear)
\startpage--\finishpage\nl
Published: \publishdate}
\vglue 0.5truein plus 0.4fil minus 0.1truein

{\parskip=0pt\leftskip 0pt plus 1fil\def\\{\par\smallskip}{\ifplaintex\large
\else\Large\fi\bf\thetitle}\par\medskip}   

\vglue 0pt plus 0.1fil 

{\parskip=0pt\leftskip 0pt plus 1fil\def\\{\par}{\sc\theauthors}
\par\medskip}

\vglue 0pt plus 0.1fil 

{\small\parskip=0pt\let\newline\\
{\leftskip 0pt plus 1fil\def\\{\par}{\sl\theaddress}\par}
\expandafter\ifx\theemail\relax    
\relax\else\vglue 5pt plus 0.02fil minus 2pt\def\\{\stdspace{\rm 
and}\stdspace} 
\cl{Email:\stdspace\tt\theemail}\fi
\ifx\theurl\relax                  
\relax\else\vglue 5pt plus 0.02fil minus 2pt\def\\{\stdspace{\rm 
and}\stdspace}
\cl{URL:\stdspace\tt\theurl}\fi\par}

\vglue 7pt plus 0.3fil minus 3pt

{\bf Abstract}
\vglue 5pt plus 0.1fil minus 2pt

\theabstract

\vglue 7pt plus 0.3fil minus 3pt

{\bf AMS Classification numbers}\quad Primary:\quad \theprimaryclass

Secondary:\quad \thesecondaryclass

\vglue 5pt plus 0.3fil minus 2pt

{\bf Keywords:}\quad \thekeywords

\vglue 10pt plus 0.5fil minus 5pt

{\small  Proposed: \theproposer\hfill Received: \receiveddate\nl
Seconded: \theseconders\hfill 
\ifx\reviseddate\relax                         
Accepted: \accepteddate                        
\else
Revised: \reviseddate                          
\fi}
\eject
}       

\let\maketitlepage\maketitlep
\let\maketitle\maketitlepage

\font\phead=cmsl9 scaled 950
\font=cmsl9 scaled 1050
\font\pnum=cmbx10 scaled 913
\font\lnum=cmbx10 
\font\pfoot=cmsl9 scaled 950
\font=cmsl9 scaled 1050
\ifplaintex
\headline{\vbox to 0pt{\vskip -4.5mm\line{\small\phead\ifnum
\count0=\startpage ISSN 1364-0380 (on line)
1465-3060 (printed) \hfill {\pnum\folio}\else\ifodd\count0\def\\{ }%
\ifx\theshorttitle\relax\thetitle\else\theshorttitle\fi\hfill{\pnum\folio}
\else\def\\{ and }{\pnum\folio}\hfill\ifx\theshortauthors\relax\theauthors
\else\theshortauthors\fi\fi\fi}\vss}}
\footline{\vbox to 0pt{\vglue 0mm\line{\small\pfoot\ifnum\count0=\startpage
\copyright\ \gtp\hfill\else
\gt, Volume \thevolumenumber\ (\thevolumeyear)\hfill\fi}\vss
}}
\else
\makeatletter
\def\@oddhead{{\small\lhead\ifnum\count0=\startpage ISSN 1364-0380 (on line)
1465-3060 (printed) \hfill {\lnum\number\count0}\else\ifodd\count0
\def\\{ }\ifx\theshorttitle\relax \thetitle \else\theshorttitle\fi\hfill
{\lnum\number\count0}\else\def\\{ and }{\lnum\number\count0}
\hfill\ifx\theshortauthors\relax 
\theauthors\else\theshortauthors\fi\fi\fi}}\def\@evenhead{\@oddhead}
\def\@oddfoot{\small\lfoot\ifnum\count0=\startpage\copyright\ \gtp\hfill\else
\gt, Volume \thevolumenumber\ (\thevolumeyear)\hfill\fi}
\def\@evenfoot{\@oddfoot}
\makeatother
\fi

\newwrite\gtoutfile
\long\gdef\makeheadfile{  
{\def\\{, }\def\s{ }
\immediate\openout\gtoutfile head.xxx
\immediate\write\gtoutfile{Proxy-for: \ifx\theasciiauthors\relax
\theauthors\else\theasciiauthors\fi\s<\ifx\theasciiemail\relax\theemail\else\theasciiemail\fi>}
\immediate\write\gtoutfile{\noexpand\\}
\immediate\write\gtoutfile{Authors: \ifx\theasciiauthors\relax
\theauthors\else\theasciiauthors\fi}
{\def\\{ }\immediate\write\gtoutfile{Title: \ifx\theasciititle\relax
\thetitle\else\theasciititle\fi}}
\immediate\write\gtoutfile{Subj-class: GT or SG or MG etc}
\immediate\write\gtoutfile{MSC-class: \theprimaryclass\ifx\thesecondaryclass\relax\else, \thesecondaryclass\fi}
\immediate\write\gtoutfile{Journal-ref: Geom. Topol. \thevolumenumber
(\thevolumeyear) \startpage-\finishpage}
\immediate\write\gtoutfile{Comments: Published by Geometry and Topology at}
\immediate\write\gtoutfile{\s\s http://www.maths.warwick.ac.uk/gt/GTVol\thevolumenumber/paper\thepapernumber.abs.html}
\immediate\write\gtoutfile{\noexpand\\}
\immediate\write\gtoutfile{}
\ifx\theasciiabstract\relax
\immediate\write\gtoutfile{\theabstract}\else
\immediate\write\gtoutfile{\theasciiabstract}\fi
\immediate\write\gtoutfile{}
\immediate\write\gtoutfile{\noexpand\\}
\immediate\write\gtoutfile{}
\immediate\closeout\gtoutfile}}  

\def\maketitlepage{\maketitlep\makeheadfile}
\let\maketitle\maketitlepage

\lognumber{473}
\received{8 July 2004}
\volumenumber{9}\papernumber{10}\volumeyear{2005}
\pagenumbers{315}{339}   
\revised{23 January 2005}
\published{23 February 2005}
\accepted{21 February 2005}
\proposed{Ronald Fintushel}
\seconded{Walter Neumann, Ronald Stern}

\usepackage{amsmath,amssymb}

\theoremstyle{plain}
\newtheorem{thm}{Theorem}[section]

\newtheorem{lem}[thm]{Lemma}

\newtheorem{conj}[thm]{Conjecture}
\newtheorem{example}[thm]{Example}
\theoremstyle{remark}
\newtheorem{rmk}[thm]{Remark}
\begin{document}
\title[Periodic maps of composite order]{Periodic maps of composite order\\ 
on positive definite 4--manifolds}
\asciititle{Periodic maps of composite order on positive definite 4-manifolds}
\author{Allan L Edmonds}
\address{Department of Mathematics, Indiana University\\Bloomington, 
IN 47405, USA}
\email{edmonds@indiana.edu}

\begin{abstract}
The possibilities for new or unusual kinds of topological, locally
linear periodic maps of \emph{non-prime} order on closed, simply
connected $4$--manifolds with positive definite intersection pairings
are explored.  On the one hand, certain permutation representations on
homology are ruled out under appropriate hypotheses.  On the other
hand, an interesting homologically nontrivial, pseudofree, action of
the cyclic group of order 25 on a connected sum of ten copies of the
complex projective plane is constructed.
\end{abstract}
\asciiabstract{%
The possibilities for new or unusual kinds of topological, locally
linear periodic maps of non-prime order on closed, simply
connected 4-manifolds with positive definite intersection pairings
are explored.  On the one hand, certain permutation representations on
homology are ruled out under appropriate hypotheses.  On the other
hand, an interesting homologically nontrivial, pseudofree, action of
the cyclic group of order 25 on a connected sum of ten copies of the
complex projective plane is constructed.}

\asciikeywords{Periodic map, 4-manifold, positive definite, permutation 
representation, pseudofree}

\maketitle
\section{Introduction}The purpose of this paper is to explore
the possibilities for unfamiliar kinds of topological, locally linear periodic
maps of \emph{non-prime} order on closed, simply connected
$4$--manifolds with positive definite intersection pairings.

Based on an earlier detailed study in Edmonds \cite{Edmonds97} of group actions
on the
$E_8$ $4$--manifold, we conjectured that a locally linear
group action on a closed, simply connected
$4$--manifold with positive definite intersection pairing induces a
(perhaps signed) permutation representation on integral homology.
For purely algebraic reasons these are exactly the kinds of
representations that occur as automorphism groups of the standard
pairing $n\left<+1\right>$.  We will in this paper, therefore,
concentrate on the problem of what kinds of permutation
representations can be realized, generally ignoring, however, the
issue of nontrivially signed permutations, by restricting to odd periods.  We will most often also
concentrate on the prototypical case of $\#_n\mathbf{C}P^2$.  We
remark that the condition of being a permutation representation
corresponds to the geometrical condition that the singular set for
the group action contains no surfaces of positive genus. See \cite{Edmonds89}, Proposition 2.4, for example.

It is well known and easy to see that every element of the
permutation group
${\Sigma}_n$ acting on $n\left<+1\right>$ can be realized by a homeomorphism (even a
diffeomorphism) of
$\#_n\mathbf{C}P^2$.  Our question in this context becomes: Which
elements of 
${\Sigma}_n$ can be realized by a periodic homeomorphism? Of
course, the answer to this question depends only on the conjugacy
class of the permutation.

We will  show that certain homology permutation representations do not arise. We also construct a new topological, locally linear, action of the cyclic group $C_{25}$ that has only a discrete singular set but is not semifree.  It is possible, however, that all three
conjectures hold as stated for smooth actions.  See Section \ref{sec:results} for precise statements of theorems.

After an earlier version of this work was written M~Tanase
\cite{Tanase2003} and I~Hambleton and Tanase
\cite{HambletonTanase2004} addressed some of the problems raised here
in the smooth category. They showed using the equivariant Yang--Mills
theory of Hambleton and R~Lee \cite{Hambleton95} that certain
permutation representations cannot be realized smoothly, including the
action of $C_{25}$ mentioned above. They showed that any such smooth
action has the same equivariant form and fixed point data of a
standard equivariant connected sum of linear actions on copies of
$\mathbf{C} P^{2}$. And they showed that an action that is smooth and
pseudofree must be semifree, in contrast to the situation in the
topological category that we display.

\section{Notation and terminology} 
We refer to \cite{Bredon72} and \cite{allday93} for generalities
about transformation groups. 

The conjugacy class of a permutation is determined by the
associated partition of
$n$ given by the cycles in the permutation.   We might write
$(p_1)( p_2) \ldots (p_r)$ where the $p_i\ge 1$ and $\sum p_i=n$. 
So, which conjugacy classes $(p_1)( p_2) \ldots (p_r)$ of
$\Sigma_n$ can be realized by a periodic map $T$ of order
$m=\text{lcm}\{p_1,\ldots, p_r\}$?

Such an action of the cyclic group $C_m$ on the set $\{1,2,\dots,n
\}$ determines a permutation representation of $C_m$ on
$\mathbf{Z}^n$. This representation we describe additively as 
$(p_1)+(p_2)+\cdots+(p_r)$

All group actions in this paper will be \emph{effective,} so that
the only group element that acts as the identity is the identity
element.  If a group $G$ acts on a space or just a set
$X$ and $x\in X$, then the
\emph{orbit} of $x$ is $G(x)$.  As a $G$--space $G(x)\approx G/H$
where
$H=G_x$, the \emph{isotropy group} of $x$.  We refer to $G/H$, as
the  ``orbit type'', and to $H$ as the corresponding ``isotropy
type''.  When $G$ is cyclic, we can abbreviate subgroups and
quotients to their integer orders. We will need to relate the
various orbit types of an action to the conjugacy types of the
corresponding representation on homology.

All group actions discussed in this paper will be locally
linear, so that each point $x$ has a neighborhood invariant under
the isotropy group $G_x$ on which the action of $G_x$ is
equivalent to a linear action on a suitable euclidean space.

If $G$ acts locally linearly on a 4--manifold $X$ with cyclic isotropy
groups of odd order, then each singular point $x\in X$ has a local
representation type, which can be described as a 
complex representation of the corresponding isotropy group.  If $x$
is a fixed point of
$G_x\approx C_m$, then this local representation can be described by
the corresponding
\emph{fixed point data}, which is an ordered pair $(a,b)$, where a
choice of generator $T\in G_x$ acts as $(z,w)\to (\zeta^a z, \zeta^b
w)$, where $\zeta=\exp(2\pi i/m)$, in local complex coordinates.

A group action is \emph{semifree} if the only isotropy group that
appears is the whole group, ie, if a nontrivial group element
fixes a point, then every group element fixes that point. A cyclic
group of prime order necessarily acts semifreely.  A group action is
said to be
\emph{pseudofree} if each nontrivial group element has a discrete
fixed point set.  Pseudofree actions are easier to study and
amenable to geometric constructions and the application of
surgery-based techniques, especially on 4--manifolds, since the
fundamental group of the complement of the singular set is the same
as the fundamental group of the ambient manifold.

\section{Basic examples}\label{sec:examples} Here we record familiar basic linear actions
of a finite cyclic group $C_m$ of order $m$ on projective spaces

Let $T$ be a periodic map of order $m$, which acts linearly on
$S^4\subset
\mathbf{C}\oplus \mathbf{C}\oplus \mathbf{R}$ by 
$T(x,y,t)=(\zeta^a x,\zeta^{b} y, t)$ or on
$\mathbf{C}P^2$ by $T[x,y,z]=[\zeta^a x,\zeta^{b} y, z]$, where
$\zeta=\exp(2\pi i/m)$.  Now $\zeta^a$ has order $m/\gcd(a,m)$ and
$\zeta^b$ has order $m/\gcd(b,m)$.  
To guarantee that $T$ actually has order $m$, and not less, we must require that $\gcd(a,b,m)=1$.

One obtains periodic maps on
$S^4$ and 
$\mathbf{C}P^2$ with at most four orbit types (in the case of $S^{4}$) and at most five orbit types (in the case of $\mathbf{C}P^2$), depending on
the choice of $m$, $a$, and $b$.  

If $a\not\equiv b\mod m$, then $T$ has exactly three fixed points: 
$[0,0,1]$, $[0,1,0]$, $[1,0,0]$ with local fixed point data $(a,b)$,
$(a-b,-b)$, and $(b-a,-a)$, respectively.

If $a\equiv b\mod m$, then $T$ fixes the point  
$[0,0,1]$ and the 2--sphere $[x,y,0]$.  The isolated point has local
fixed point data
$(a,a)$, the 2--sphere has normal euler number $+1$, and the
generator $T$ rotates by
$2\pi a/m$ in the normal fiber.

Equivariant (anti-holomorphic) blow up of orbits
yields similar actions with arbitrarily many singular points on
$\#_n\mathbf{C}P^2$.  One can write down the fixed point data for
such a positive blow-up as follows.  Start with the data
$$(a,b),(a-b,-b),(b-a,-a).$$
Using the negative of the second entry we also have $\mathbf{C}P^2$
data
$$(a-b,b),(a-2b,-b),(2b-a,b-a).$$
Equivariant connected sum yields an action on $\#_2\mathbf{C}P^2$
with data
$$(a,b),(b-a,-a),(a-2b,-b),(2b-a,b-a).$$
One can then repeat the process. At each stage there are several options for the choice of point to blow up. Tanase \cite{Tanase2003} and Hambleton and Tanase \cite{HambletonTanase2004} have given a thorough analysis of the possibilities for the fixed point data and the permutation representations that can arise when one iterates this construction.

\section{Statements of results}\label{sec:results} We give three results here that show that under appropriate hypotheses certain permutation representations do not arise from topological, locally linear actions on a connected sum of copies of $\mathbf{C}P^{2}$.

\begin{thm}\label{thm:best}For $p$ an odd
prime, there is no locally linear action of C$_m$,
$m=p^k$, $k\ge 3$, on $\#_n\mathbf{C}P^2$ (or any closed, simply
connected, positive definite 4--manifold with $b_{2}=n$), , where $n=p^{2}+p$ and the representation on homology is of type
$(p)+(p^2)$.
\end{thm}

\begin{thm}\label{thm:p^k}For $p$ an odd prime, there is
no effective, locally linear, pseudofree, action of C$_m$,
$m=p^k$, $k\ge 2$, on $\#_n\mathbf{C}P^2$ (or any closed, simply
connected, positive definite 4--manifold) with homology representation
of type
$$(p^\ell)+r_{\ell+1}(p^{\ell+1})+\dots+r_k(p^k)$$ where
$1\le\ell\le k-1$.  
\end{thm}

 The proof involves the careful study of the $G$--Signature
Formulas for $T$ and
$T^p$, where $T$ denotes a chosen generator of $C_m$. 

\begin{thm}\label{thm:p=3} There is no locally linear,
\emph{pseudofree} action of C$_9$ on
$\#_n\mathbf{C}P^2$ (or any closed, simply
connected, positive definite 4--manifold) with homology
representation of type
$r(9)+s(3)+t(1)$, where
$t\le 1$ and $s\ge 1$.
\end{thm} 
When $t=0$ and $s=1$, this is a consequence of Theorem \ref{thm:p^k},
but it requires more work when $t=1$ or $s\ge 2$. The proof again
involves a detailed look at the $G$--Signature Formulas for
$T$ and
$T^3$, where
$T$ generates
$C_9$.  But the argument is more subtle.

In contrast to these non-existence results, we have one successful
realization theorem.

\begin{thm}\label{thm:p=5} There is a locally linear {pseudofree}
action of C$_{25}$ on
$\#_{10}\mathbf{C}P^2$ with homology representation of type
$2(5)$.
\end{thm} 
This result is especially interesting in light of the work of Tanase \cite{Tanase2003} and of Hambleton and Tanase \cite{HambletonTanase2004} which shows that such an action cannot be smoothed.
\begin{rmk}
From this result it is easy to use equivariant blow up to construct actions of C$_{25}$ on
$\#_{n}\mathbf{C}P^2$ with homology representation of type
$r(25)+2(5)+t(1)$, for any
$r\ge 0$ and $t\ge 0$, where $n=10+25r+t$.
\end{rmk}

In this case, suitable data that might come from an actual action
is produced, in the sense that the $G$--Signature Formula would be
satisfied. Then techniques developed in Edmonds and Ewing \cite{EdmondsEwing92} are
used to show that the data can actually be realized by a group
action.  One common thread of all these  results is that in some
sense the $G$--Signature Formula holds most of the key to the
existence of group actions.

\section{Tools} Because they are crucial in what follows, we explicitly record two important and well-known tools
for use in subsequent sections.

\subsection{Lefschetz Fixed Point Formula}

Let $T\co X\to X$ generate an action of $C_m$ on $X$, a closed, simply
connected 4--manifold.  Then ``local Smith theory'' implies that
$F=\text{Fix}(T)$ consists of isolated points and surfaces.  And the
Lefschetz fixed point formula says that the Euler characteristic of the fixed point set is given by
$$
\chi(F) =\Lambda(T)=2+\text{trace}\,[\,T_*\co H_2(X)\to H_2(X)\,].
$$
For a
general discussion of this basic result see
Allday and Puppe \cite{allday93}, (3.29). A refinement using a closer study of
equivariant cohomology shows that all fixed surfaces are 2--spheres
if and only if the representation on $H_2$ is a permutation
representation.  See Edmonds \cite{Edmonds89}, Proposition 2.4.

\subsection{$G$--Signature Formula} Let $T\co X\to X$ generate an action
of $C_m$ on a  4--manifold $X$.  In general, we understand the ``$g$--signature'' to be the character of a certain virtual complex representation evaluated on the generator $T$.  One extends the intersection pairing on $H_2(X)$ to a Hermitian pairing on $H_2(X;\mathbf{C})$ and forms the differences of the traces of $T_*$ on the positive and negative parts.    But in the case of interest in this paper of positive definite $4$--manifolds, it simplifies dramatically to
$$
\sigma(T,X) := \text{trace}\,[\,T_*\co H_2(X)\to H_2(X)\,]
$$ 
and assuming the representation on $H_2(X)$ is of permutation type, it is the number of fixed basis vectors in the permutation representation.  In particular, it is a rational integer in this case.

Set $\zeta=\text{exp}(2\pi i/m)$.  Suppose $T$ has isolated fixed
points $x_i$ and fixed surfaces $S_j$.  Suppose $T$ has local
representation of ``type'' $(a_i,b_i)$ at $x_i$; let $S_j$ have normal
euler number $n_j$ and normal rotation angle data $e_j$ (so that $T$
rotates an oriented normal plane to $S_j$ by $2\pi e_j/m)$).  Then we
have:
\begin{gather*}
\sigma(T,X) =
\sum_{i}^{}\frac{(\zeta^{a_i}+1)}{(\zeta^{a_i}-1)}\frac{(\zeta^{b_i}+1)}
{(\zeta^{b_i}-1)} -
\sum_j\frac{4n_j\zeta^{e_j}}{(\zeta^{e_j}-1)^2} 
\end{gather*}
A nice reference for this version is Gordon \cite{Gordon86}.

\section{Number-theoretic issues} This section may be skimmed or
skipped and referred back to as needed.  The issue is manipulating
and drawing conclusions from the $G$--Signature Formula. The key results are Theorems \ref{thm:t=0} and \ref{thm:2}, which are at least on the surface reminiscent of the Franz Independence Lemma \cite{Franz35}.

 We begin with a few very
elementary but useful lemmas.  If
$\zeta=\exp(2\pi i/m)$, then  we can identify the
cyclotomic field $\mathbf{Q}(\zeta)$ with
$\mathbf{Q}[x]/\Phi_m(x)$, where $\Phi_m(x)$ is the $m$th
cyclotomic polynomial.  Any elementary number theory book, for example \cite{BorevichShafarevich66}, should be
adequate reference.

\begin{lem} If $\zeta=\exp(2\pi i/m)$, then $\zeta^a-1$
($a\not\equiv 0\mod m$) has inverse in
$\mathbf{Q}(\zeta)$ given by
$$
(\zeta^a-1)^{-1}=\frac{-1}{m}\prod_{\substack{i=1\\ (i\ne
a)}}^{m-1}(1-\zeta^i).
$$
\end{lem}
\begin{proof}
We have the complex polynomial factorization
$$
(x^m-1) = \prod_{i=0}^{m-1}(x-\zeta^i).
$$
Dividing through by $x-1$ we have
$$
1+x+\cdots + x^{m-1}= \prod_{i=1}^{m-1}(x-\zeta^i).
$$
Plug in $x=1$ to get
$$
m=\prod_{i=1}^{m-1}(1-\zeta^i)
$$
and isolate $\zeta^a-1=-(1-\zeta^a)$ to get the result.
\end{proof}

\begin{lem}
If
$$
\frac{(\zeta^{a}+1)}{(\zeta^{a}-1)} = \frac{(\zeta^{b}+1)}
{(\zeta^{b}-1)}
$$
where $\zeta = \exp(2\pi i/m)$, then $a\equiv b \mod m$.
\end{lem}
Cross-multiply and simplify to get $\zeta^a=\zeta^b$.  The
result follows.

\begin{lem}
If $\zeta = \exp(2\pi i/m)$, then:
$$
\frac{(\zeta^{-a}+1)}{(\zeta^{-a}-1)} =
-\frac{(\zeta^{a}+1)}{(\zeta^{a}-1)}
$$
\end{lem}
\begin{proof}
Multiply numerator and denominator on the left by $\zeta^a$ and
simplify.
\end{proof}

\begin{lem}\label{lem:1}
The equation
$$
\frac{(\zeta^{a}+1)}{(\zeta^{a}-1)}\frac{(\zeta^{b}+1)}
{(\zeta^{b}-1)} = 1
$$
in which $\zeta=\exp(2\pi i$m), $m$ odd,
has no solutions $a$ and $b$. 
\end{lem}

\begin{proof}
Clear denominators and simplify to show that $\zeta^b=-\zeta^a$.
The result follows.
\end{proof}

\begin{thm}\label{thm:t=0} Let $m=p^t$, where $p$ is an
odd prime. If $\zeta=\exp(2\pi i/m)$ and
$$
\frac{(\zeta^{a}+1)}{(\zeta^{a}-1)}
\frac{(\zeta^{b}+1)}{(\zeta^{b}-1)} =
\frac{(\zeta^{c}+1)}{(\zeta^{c}-1)}
\frac{(\zeta^{d}+1)}{(\zeta^{d}-1)}
$$ in $\mathbf{Q}(\zeta)$, then
$\{c,d\}\equiv\pm\{a,b\}\mod m$.
\end{thm}

\begin{rmk} It seems doubtful that Theorem \ref{thm:t=0} is  true
for any odd integer $m$.  But the proof would definitely be more
difficult in that case.  Note that by applying a Galois automorphism of $\mathbf{Q}(\zeta)/\mathbf{Q}$, the hypothesis implies that the formula holds for $\zeta$ replaced by any \emph{primitive} $m$th root of unity.  As a corollary, the formula holds for $\zeta$ replaced by any $m$th root of unity at all.
\end{rmk}

\begin{rmk} Theorem \ref{thm:t=0} is equivalent to the
statement that if
$$
\frac{(\zeta^{a}+1)}{(\zeta^{a}-1)}
\frac{(\zeta^{b}+1)}{(\zeta^{b}-1)} +
\frac{(\zeta^{c}+1)}{(\zeta^{c}-1)}
\frac{(\zeta^{d}+1)}{(\zeta^{d}-1)} =0
$$ then $\{c,d\}\equiv\pm\{-a,b\}\mod m$, which is the form in which we shall apply it.
\end{rmk}

\begin{proof} If $a=\pm c$ or $a=\pm d$ or $b=\pm c$ or
$b=\pm d$, then the result follows from our
analysis of simpler cases that arise when we divide both sides by
common factors.

Cross-multiply  and simplify to
obtain:
\begin{center}
\begin{tabular}{|c|}
\hline
 \\
$
\ \zeta^{a+c+d}+\zeta^{b+c+d}+\zeta^a+\zeta^b =
\zeta^{a+b+c}+\zeta^{a+b+d}+\zeta^c+\zeta^d\ 
$
\\
\\
\hline
\end{tabular}
\end{center}
It will suffice to show that some exponent on the left equals some exponent on the right (mod $m$).  To see this, focus, for example, on
$a$.  If it equals $a+b+c \mod m$, then $b\equiv
-c$ and we are done, by the preceding observation.  If it equals $a+b+d$, then $b\equiv
-d$, and again we are done.  If it equals $c$, then
$a\equiv c$, and again we are done.
And, finally, if  it equals $d$, then $a\equiv d$, and again we are done.  Similar considerations apply to any of the other exponents on the left.  Altogether there would be 16 similar cases, which can be checked in ones head.

Now any element of $\mathbf{Q}(\zeta)$ has a unique
expression as a polynomial in $\zeta$ of degree less than
$\deg(\Phi_m(x))$.  In the special case when all the
indicated exponents are less than $\deg(\Phi_m(x))$, the
result follows, since all the exponents on the left side must
coincide with exponents on the right side.  

The subsequent argument aims at showing that some exponent on the
left must equal some exponent mod $p^t$ on the right in all cases.

Fix $m=p^t$.  Let us say that an integer $k$ is
\emph{ordinary}
mod $m$ if its minimal non-negative representative
$\overline{k}$ mod $m$ is less than
$\deg(\Phi_m(x))=p^t-p^{t-1}=(p-1)p^{t-1}$.  In this
case the unique expression for $\zeta^k$ is simply
$\zeta^{\overline{k}}$.

Similarly, let us say  than an integer $k$
is \emph{critical} mod $m$ if its minimal non-negative
representative
$\overline{k}$ is greater than or
equal to $\deg(\Phi_m(x))$, ie, $(p-1)p^{t-1}\le
\overline{k}\le p^t-1$ We must describe the unique
representative of
$\zeta^k$.

Recall that for $m=p^t$ the cyclotomic polynomial
$$\Phi_m(x)=1+x^{p^{t-1}}+x^{2{p^{t-1}}}+\dots
+x^{(p-1){p^{t-1}}}.$$
If $(p-1)p^{t-1}\le
\overline{k}\le p^t-1$, then
$\overline{k}=(p-1)p^{t-1}+\ell$, where, of course,
$\ell=\overline{k}-(p-1)p^{t-1}$. Then the relation
$x^\ell\Phi_m(x)=0$, holds in $\mathbf{Q}[x]/\Phi_m(x)$,
that is
$$
x^\ell+x^{p^{t-1}+\ell}+x^{2p^{t-1}+\ell}+\dots
+x^{(p-1)p^{t-1}+\ell}=0
$$
or
$$
x^k=-x^{k-p^t+p^{t-1}}-x^{k-p^t+2p^{t-1}}-\dots
-x^{k-p^t+(p-1)p^{t-1}}.
$$
Thus when $k$ is critical mod $m$, the unique
representative of $\zeta^k$ contributes $p-1$ terms,
each with coefficient $-1$.  Note also that if $k_1$ and
$k_2$ are both critical, but distinct mod $m$, then the
unique representatives of $\zeta^{k_1}$ and
$\zeta^{k_2}$ have no terms in common.  This crucial
observation is not at all clear in the case when $m$ is
not a prime power.

Now let us return to our equation in the form
$$
\zeta^{a+c+d}+\zeta^{b+c+d}+\zeta^a+\zeta^b =
\zeta^{a+b+c}+\zeta^{a+b+d}+\zeta^c+\zeta^d.
$$
Since we already discussed the case when
all exponents on the left side are ordinary, we can
assume there is at least one critical exponent on the
left hand side.  This means that in the unique
representation for the left hand side there
is at least one negative term (e.g., there are
actually two or three critical terms, or just one
but $p-4\ge 1$ ).  Then the right hand side also has at
least one negative term.  But from the knowledge of the
degree of a negative term we can determine the original
degree (since $m$ is a prime power).  Thus at least one
exponent on the right hand side equals an exponent on
the left hand side. Then we are done, as we observed at
the beginning.

Strictly speaking we have to look more closely when
$p=3$, where the left hand side has exactly one
critical exponent and three ordinary terms exactly
cancel off the $p-1=2$ negative terms from the
representation of the critical term.  But then it
follows that the right hand side has exactly one
positive ordinary term left, just as the left hand
side must have.  In particular there is an exponent that
appears on both sides of the original equation, as needed.
\end{proof}

Similar considerations are used to prove the following result.

\begin{thm}\label{thm:2} Let $m=p^t$, where $p$ is an
odd prime.  Suppose that $k\in\mathbf{Z}$ and that $c$
and $d$ are integers prime to $p$. Set
$\zeta=\exp(2\pi i/m)$.  If
$$
\frac{-4k\zeta}{(\zeta-1)^2}
+p\frac{(\zeta^{c}+1)}{(\zeta^{c}-1)}
\frac{(\zeta^{d}+1)}{(\zeta^{d}-1)}
=
p
$$ 
in $\mathbf{Q}(\zeta)$, then
$k=p$, and $c\equiv d\equiv 1\mod p^t$ or $c\equiv d\equiv -1\mod p^t$.
\end{thm}
\begin{rmk}
When $k=p$, and $c\equiv d\equiv 1\mod p^t$, the formula is just $p$
times the $G$--Signature Formula for a standard periodic map on the
complex projective plane.  The theorem shows that there are no
mysterious solutions to the $G$--signature equation in this case.  It then follows that
$$
\frac{-4k\xi}{(\xi-1)^2}
+p\frac{(\xi^{c}+1)}{(\xi^{c}-1)}
\frac{(\xi^{d}+1)}{(\xi^{d}-1)}
=
p
$$ 
for any $\xi=\zeta^{p^j}$, $j<t$.

\end{rmk}

\begin{proof}We may assume that $k>0$ by Lemma \ref{lem:1}.  By replacing $(c,d)$ by $(d,c)$ and/or $(p^t-c,p^t-d)$, we may assume that $1\le c\le d\le p^t-1$ and $c\le \frac{p^t-1}{2}$. 

Observe that $k$ must be divisible by $p$, as we can see by multiplying through by ${(\zeta^{c}-1)}{(\zeta^{d}-1)}$ and reducing modulo ${(\zeta-1)}$.
Therefore one may write $k=p\ell$.  Dividing through by $p$, then, we must
solve the equation
$$
\frac{-4\ell\zeta}{(\zeta-1)^2}
+\frac{(\zeta^{c}+1)}{(\zeta^{c}-1)}
\frac{(\zeta^{d}+1)}{(\zeta^{d}-1)}
=1.
$$
Clearing denominators, multiplying some things out, cancelling some
terms and collecting and factoring, we obtain
$$
2\ell\zeta(\zeta^c-1)(\zeta^d-1) = (\zeta -
1)^2(\zeta^{c}+\zeta^{d}).
$$
(As a quick check note that $c=d=\ell=1$ does satisfy this equation
and moreover, if $c=d=1$, then it follows that $\ell=1$.)

Now multiply it all out and get:
\begin{center}
\begin{tabular}{|c|}
\hline
 \\
\ $2\ell\zeta^{c+d}+2\ell=\zeta^{c+1}+\zeta^{d+1}+(2\ell-2)\zeta^{c}+(2\ell-2)\zeta^{d}+\zeta^{c-1}+\zeta^{d-1}$\ 
\\
\\
\hline
\end{tabular}
\end{center}
When the left hand side and right hand side are each put into their unique form as a polynomial of degree less than $p^t-p^{t-1}$, we can compare constant terms.  On the left hand side we either have $2\ell\ge 2$ or $0$, and the latter only if $c+d\equiv -p^{t-1} \mod p^t$.  

Suppose $2\ell\ge 2$ is the constant term on the left hand side. On the right hand side, then, the constant term can only have a positive contribution from a 0 power of $\zeta$.    This could occur only from $c\equiv 1$ or $c\equiv-1$ or $d\equiv1$ or $d\equiv-1$ (mod $p^t$).  This would limit us to $\ell=1$ and two of these four cases occuring.  If $c\equiv d\equiv\pm 1 \mod p^t$, we are done.  The other possibility is $c\equiv 1\mod p^t$ and $d\equiv p^t-1\mod p^t$.  But then $c+d$ becomes critical and the constant term on the left hand side would be 0.  

Thus we finally have to consider the case when the constant term on the left hand side is 0, and somehow eliminate this possibility.  That happens only when $c+d\equiv p^t-p^{t-1}\mod p^t$.  This clearly leads to a quick contradiction, but it remains to be seen how to organize the proof efficiently.  In this situation the left hand side would have exactly $p-2$ terms in its reduced form:
$$
-2\ell(\zeta^{p^{t-1}}+\zeta^{2p^{t-1}}+\dots
+\zeta^{(p-2)p^{t-1}})
$$

Since $c\ge 1$ and $c+d=p^t-p^{t-1}$, we see that 
$$d=p^t-p^{t-1}-c\le p^t-p^{t-1}-1$$
or
$$d+1=p^t-p^{t-1}-c+1\le p^t-p^{t-1}.$$
Thus the only possible critical exponent on the right hand side is $d+1$. And $d+1$ is critical only if 
$$d+1=p^t-p^{t-1}-c+1= p^t-p^{t-1}$$
if and only if
$c=1$.  Thus, if there is a critical exponent on the right hand side, the right hand side can be rewritten as
$$
\zeta^{2}-1-\zeta^{p^{t-1}}\!\!-\zeta^{2p^{t-1}}\!\!-\dots-
\zeta^{(p-2)p^{t-1}}\!
+(2\ell-2)\zeta^{1}+(2\ell-2)\zeta^{p^t-p^{t-1}-1}+\zeta^{0}+\zeta^{p^t-p^{t-1}-2}
$$
or
$$
(2\ell-2)\zeta^{1}+\zeta^{2}-\zeta^{p^{t-1}}\!-\zeta^{2p^{t-1}}\!-\dots-\zeta^{(p-2)p^{t-1}}
+(2\ell-2)\zeta^{p^t-p^{t-1}-1}+\zeta^{p^t-p^{t-1}-2}.
$$
Thus on the right hand side the coefficient of $\zeta$ is $2\ell-2$, while on the left hand side there is no $\zeta$ term unless $t=1$.  But then the coefficient of $\zeta$ on the left hand side is $-2\ell$.  We conclude that $2\ell-2=-2\ell$, or $\ell={1}/{2}$, a contradiction.

Finally we have to consider the case where no exponent on the right hand side is critical.  
$$
-2\ell(\zeta^{p^{t-1}}\!\!+\zeta^{2p^{t-1}}\!\!+\dots
+\zeta^{(p-2)p^{t-1}})=\zeta^{c+1}\!+\zeta^{d+1}\!+(2\ell-2)\zeta^{c}\!+(2\ell-2)\zeta^{d}\!+\zeta^{c-1}\!+\zeta^{d-1}
$$
In this case the right hand side has at least 3 distinct terms and at most 6.  But the left hand side has $p-2$ terms.  Thus $3\le p-2\le 6$.  Meanwhile on the left hand side all terms have the same coefficient of $-2\ell$.  The only way for some terms on the right hand side to coalesce would be to have $c=d$.  (The other possibility of $c=d-1$ is ruled out because $c+d=p^t-p^{t-1}$ is even.)  So our equation finally becomes
$$
-2\ell(\zeta^{p^{t-1}}+\zeta^{2p^{t-1}}+\dots
+\zeta^{(p-2)p^{t-1}})=2\zeta^{c+1}+(4\ell-4)\zeta^{c}+2\zeta^{c-1}
$$
where $c=(p^t-p^{t-1})/{2}$ and $p=5$.  But then $-2\ell=2$ and $-2\ell=4\ell-4$, a contradiction. 
\end{proof}

When manipulating the
$G$--Signature Formula it once seemed natural to invoke the following
``obvious fact''.
\begin{conj}[$m,p,t$]\label{conj_t}
$$\displaystyle
\sum_{i=1}^{t+2}\frac{(\zeta^{a_i}+1)}{(\zeta^{a_i}-1)}\frac{(\zeta^{b_i}+1)}
{(\zeta^{b_i}-1)}=t
\leqno{\rm If}$$ where $\zeta=\exp(2\pi i/m)$, and $a_i$ and
$b_i$ are prime to $m$, then
$$\displaystyle\sum_{i=1}^{t+2}\frac{(\zeta^{pa_i}+1)}{(\zeta^{pa_i}-1)}\frac{(\zeta^{pb_i}+1)}
{(\zeta^{pb_i}-1)}=t$$ where $p|m$. 
\end{conj}

When $t=0$ this is equivalent to  Theorem
\ref{thm:t=0}. The analogous statement is easily seen to be true for
$p$ prime to
$m$, because the two formulas are related by a Galois automorphism
in that case.  We will prove this for the smallest case not covered
by earlier work:
$p=3$, 
$m=9$ and
$t=1$.  But in general we will show that this natural conjecture is
not true, for, say, $p=3$, 
$m=9$ and
$t\ge 4$.  (Perhaps it does hold for general $p$ when $t=1$, though.)

\begin{rmk}
The special case above when $t=1$ corresponds to the
situation of cyclic groups acting on $\mathbf{C}P^2$.  In that case,
under the assumption that the given equation holds for all $m$th
roots of unity, work of Edmonds and Ewing \cite{EdmondsEwing89} shows that the
\emph{only} solutions come from standard linear actions.
\end{rmk}

 \begin{thm}\label{thm:391}If
$$\displaystyle
\sum_{i=1}^{3}\frac{(\zeta^{a_i}+1)}{(\zeta^{a_i}-1)}\frac{(\zeta^{b_i}+1)}
{(\zeta^{b_i}-1)}=1$$ where $\zeta=\exp(2\pi i/9)$, and $a_i$ and
$b_i$ are prime to $9$, then
$$\displaystyle\sum_{i=1}^{3}\frac{(\zeta^{3a_i}+1)}{(\zeta^{3a_i}-1)}\frac{(\zeta^{3b_i}+1)}
{(\zeta^{3b_i}-1)}=1$$ 
\end{thm}
\begin{proof}
One can simply check this by brute force on a computer.
\end{proof}

\begin{example}\label{example:394}There exist integers (mod
9), $a_i,b_i\ \ (i=1,\ldots,8)$, such that
$$\displaystyle
\sum_{i=1}^{8}\frac{(\zeta^{a_i}+1)}{(\zeta^{a_i}-1)}\frac{(\zeta^{b_i}+1)}
{(\zeta^{b_i}-1)}=6$$ 
where $\zeta=\exp(2\pi i/9)$, and $a_i$ and
$b_i$ are prime to $9$, but
$$\displaystyle\sum_{i=1}^{8}\frac{(\zeta^{3a_i}+1)}{(\zeta^{3a_i}-1)}\frac{(\zeta^{3b_i}+1)}
{(\zeta^{3b_i}-1)}\ne 6.$$ 
\end{example}
\begin{proof}
This is not so easy to find by brute force.  But one can
reduce the problem to solving integral equations and
eventually find examples.  Of course, once one finds the
examples, the verification that the examples work is
essentially trivial.

Here is one solution:

\begin{multline*}
2\frac{(\zeta^{1}+1)}{(\zeta^{1}-1)}
\frac{(\zeta^{8}+1)}{(\zeta^{8}-1)}+
3\frac{(\zeta^{1}+1)}{(\zeta^{1}-1)}\frac{(\zeta^{2}+1)}
{(\zeta^{2}-1)}\\
+\frac{(\zeta^{1}+1)}{(\zeta^{1}-1)}\frac{(\zeta^{4}+1)}
{(\zeta^{4}-1)}+
\frac{(\zeta^{2}+1)}{(\zeta^{2}-1)}\frac{(\zeta^{7}+1)}
{(\zeta^{7}-1)}+
\frac{(\zeta^{2}+1)}{(\zeta^{2}-1)}\frac{(\zeta^{4}+1)}
{(\zeta^{4}-1)}=6
\end{multline*}

\noindent But the same formula with $\zeta$ replaced by $\zeta^3$
yields 2, not 6.
\end{proof}

\section{A non-existence result}
Here we give what in many ways is our best result, ruling out a
particular homology representation without making any extra
assumptions, such as the pseudofree hypothesis explored in the next
section,  about an action. In the smooth category the analog of this result follows from Hambleton and Tanase \cite{HambletonTanase2004}.

\begin{thm}[Restatement of Theorem \ref{thm:best}]For $p$ an odd
prime, there is no locally linear action of C$_m$,
$m=p^3$, on $\#_n\mathbf{C}P^2$, where $n=p^{2}+p$ and the representation on homology is of type
$(p)+(p^2)$.
\end{thm}

\begin{proof}
Suppose such an action exists. Let $T$ denote a
generator of
$C_m$.  The homology representation implies that $\sigma(T,X)=0$, 
$\sigma(T^{p},X)=p$ and $\sigma(T^{p^2},X)=p^2+p$.

We have inclusions
$$\text{Fix}(T)\subset\text{Fix}(T^{p})\subset\text{Fix}(T^{p^2})
$$
and $\text{Fix}(T)=\text{Fix}(T,\text{Fix}(T^{p^2}))$.
Moreover, considering the Lefschetz Fixed Point formula,
 $\text{Fix}(T)$ is two points or a 2-sphere.
Furthermore, in this case $\text{Fix}(T)$ consists either of
two fixed points (ie, $S^0)$, with cancelling local fixed point data,
or of a single 2-sphere
$S^2$ with trivial normal bundle, as we see by applying the
Lefschetz Fixed Point Formula and the $G$--Signature Formula to
$T$.

Points of $\text{Fix}(T^{p}) - \text{Fix}(T) $ come in
nontrivial orbits of size $p$.  Let $\overline{\text{Fix}(T)}$ denote the union of the components of $\text{Fix}(T^{p})$ meeting $\text{Fix}(T)$, as it is a kind of closure.  There are three cases:
\begin{enumerate}
\item $\text{Fix}(T)=S^2=\overline{\text{Fix}(T)}$
\item $\text{Fix}(T)=S^0=\overline{\text{Fix}(T)}$
\item $\text{Fix}(T)=S^0\subset \overline{\text{Fix}(T)}=S^2$
\end{enumerate}
 We will argue that
$\text{Fix}(T)=S^0$ while $\overline{\text{Fix}(T)}=S^2$ by eliminating the first two cases.  

In any case we must have that 
$$\text{Fix}(T^p)-\overline{\text{Fix}(T)}=p\ \text{points}.$$
For if
$$\text{Fix}(T^p)-\overline{\text{Fix}(T)}=aS^2+b\text{ points}$$
then the Lefschetz formula implies that $2a+b=p$.  But $T$ acts on this set without fixed points.  It follows that either $b=0$ or $b\ge p$.  But if $b=0$, then $2a=p$, contradicting the assumption that $p$ is odd.  Thus $b\ge p$, from which it follows that $b=p$ and $a=0$.  Similarly
$$\text{Fix}(T^{p^2})-\text{Fix}(T^p)=p^2\ \text{points}.$$
Suppose now that Case 1 occurs, that is, ${\text{Fix}(T)}\cong S^2$.  Then, clearly,
$\overline{\text{Fix}(T)}\cong S^2$ as well.  Now 
$$
\sigma(T,X)=0=
\frac{-4n_j\zeta^{e_j}}{(\zeta^{e_j}-1)^2}
$$
which forces the
normal euler number to be 0.  The fixed point data at each of
the $p$ new fixed points of
$T^{p}$ must be the same, since $T$ permutes them transitively.  Let it be
$(e,f)$  (mod $p^2$). Then the $G$--Signature Formula for
$T^p$ yields
$$\sigma(T^{p},X)=p=p\frac{(\zeta^{e}+1)}{(\zeta^{e}-1)}
\frac{(\zeta^{f}+1)}{(\zeta^{f}-1)}$$
which makes
$$1=\frac{(\zeta^{e}+1)}{(\zeta^{e}-1)}
\frac{(\zeta^{f}+1)}{(\zeta^{f}-1)}$$
contradicting Lemma \ref{lem:1}. We conclude that ${\text{Fix}(T)}\cong S^0$.

Suppose next that ${\text{Fix}(T)}\cong S^0$ and that
$\overline{\text{Fix}(T)}\cong S^0$ also.  Let $(a,b)$ and $(c,d)$
be the local fixed point data for $T$ (and hence $T^{p}$, too)
at the two points of $\text{Fix}(T)$.  Now the $G$--Signature Formula
for $T$ says that
$$
\frac{(\zeta^{a}+1)}{(\zeta^{a}-1)}\frac{(\zeta^{b}+1)}
{(\zeta^{b}-1)}+\frac{(\zeta^{c}+1)}{(\zeta^{c}-1)}\frac{(\zeta^{d}+1)}
{(\zeta^{d}-1)}=0
$$
where $\zeta=\exp(2\pi i/m)$.  By Theorem \ref{thm:t=0}, it follows
that the fixed point data cancels and we may assume that
$(c,d)=(a,-b)$ (mod $m$), for example.  

Again let the fixed point data at
the $p$ new fixed points of
$T^{p}$ be
$(e,f)$.  Then the $G$--Signature Formula for $T^{p}$ says that
$$
p\frac{(\xi^{e}+1)}{(\xi^{e}-1)}\frac{(\xi^{f}+1)}{(\xi^{f}-1)}
=p$$
where $\xi=\zeta^p=\exp(2\pi i/p^2)$, which cannot hold, by Lemma \ref{lem:1}.  We
therefore conclude that $\text{Fix}(T)=S^0$ while
$\overline{\text{Fix}(T)}=S^2$.

Suppose that the 2-sphere $\overline{\text{Fix}(T)}$ has normal euler number $k$. 
It then follows that the fixed point data for the two fixed points
($\text{Fix}(T)$) must have the form $(a,b)$ and $(-a,b-ka)$ (mod $m$), as
one sees by considering a standard model of the euler class $k$
bundle over $S^2$.  Compare Lemma \ref{lem:stddata} below.  On the
other hand the $G$--Signature Formula for
$T$ states that
$$
\frac{(\zeta^{a}+1)}{(\zeta^{a}-1)}\frac{(\zeta^{b}+1)}
{(\zeta^{b}-1)}+\frac{(\zeta^{-a}+1)}{(\zeta^{-a}-1)}\frac{(\zeta^{b-ka}+1)}
{(\zeta^{b-ka}-1)}=0
$$
where $\zeta=\exp(2\pi i/m)$. By Theorem \ref{thm:t=0} the data cancels.  Therefore we can
conclude that $k\equiv 0$ mod $m$.  

Then, as above, if the fixed
point data at (each of) the $p$ new fixed points of
$T^{p}$ is
$(e,f)$.  Then the $G$--Signature Formula for $T^{p}$ says that
$$
\sigma(T^p,X)=
\frac{-4k\xi^{b}}{(\xi^{b}-1)^2} +{p}\frac{(\xi^{e}+1)}{(\xi^{e}-1)}\frac{(\xi^{f}+1)}{(\xi^{f}-1)}
={p}
$$
where $\xi=\zeta^p=\exp(2\pi i/p^2)$.
This certainly has solutions, for example, when $k=p$ and $b=e=f$.  According to Theorem \ref{thm:2} these are the only solutions, up to equivalence (replacing $(e,f)$ by $(-e,-f)$ or $(f,e)$) .

Now consider $\text{Fix}(T^{p^2})$.  This consists of $\text{Fix}(T^{p})$ plus $p^2$ points making up a single $T$ orbit and $p$ orbits of $T^p$.  Let the $p^2$ new fixed points have local data $(c,d)$.  Then the g--Signature Formula for $T^{p^2}$  says that
$$
\sigma(T^{p^2},X)=
\frac{-4k\mu^{b}}{(\mu^{b}-1)^2} +{p}\frac{(\mu^{b}+1)}{(\mu^{b}-1)}\frac{(\mu^{b}+1)}{(\mu^{b}-1)}
+{p^2}\frac{(\mu^{c}+1)}{(\mu^{d}-1)}\frac{(\mu^{c}+1)}{(\mu^{d}-1)}
={p}+{p^2}$$
where $\mu=\xi^p=\zeta^{p^2}=\exp(2\pi i/p).$
It follows that
$$
{p}
+{p^2}\frac{(\mu^{c}+1)}{(\mu^{d}-1)}\frac{(\mu^{c}+1)}{(\mu^{d}-1)}
={p}+{p^2}$$
or
$$
\frac{(\mu^{c}+1)}{(\mu^{d}-1)}\frac{(\mu^{c}+1)}{(\mu^{d}-1)}
=1
$$
contradicting Lemma \ref{lem:1}.
\end{proof}

\section{Pseudofree actions} 
Here we investigate the possibilities
for a pseudofree action that is not semifree.  
We first present a result that is analogous to
the result that a prime (power) order map cannot have just one fixed
point:  a periodic map of odd prime power order cannot have just one
singular orbit of smallest size.

Recent work of Hambleton and Tanase \cite{HambletonTanase2004} shows that a pseudofree action must be semifree in the smooth category. So the actions considered here cannot be smoothed.

\begin{thm}[Restatement of Theorem \ref{thm:p^k}]
For $p$ an odd prime, there is
no effective, locally linear, pseudofree, action of C$_m$,
$m=p^k$, $k\ge 2$, on $\#_n\mathbf{C}P^2$ (or any closed, simply
connected, positive definite 4--manifold) with homology representation
of type
$$(p^\ell)+r_{\ell+1}(p^{\ell+1})+\dots+r_k(p^k)$$ where
$1\le\ell\le k-1$.  
\end{thm}

\begin{rmk}
An action with such a homology representation can exist for
non-pseudofree actions.  For example $C_{p^2}$ acts on $S^4$ with
two fixed points and a $2$--sphere of points with isotropy group
$C_p$.  Blowing up one orbit of type $C_p$ and an arbitrary
number of principal orbits creates the homology representation
$(p)+r_2(p^2)$.
\end{rmk}

\begin{proof}Suppose such an action exists. Let $T$ denote a
generator of
$C_m$.  The homology representation implies that $\sigma(T,X)=0$ and 
$\sigma(T^{p^\ell},X)=p^\ell$.

We have inclusions
$$\text{Fix}(T)=\text{Fix}(T^{p})=\cdots=\text{Fix}(T^{p^{\ell-1}})
\subset
\text{Fix}(T^{p^\ell})$$
and $\text{Fix}(T)=\text{Fix}(T,\text{Fix}(T^{p^\ell}))$.
Moreover, 
 $\text{Fix}(T^{p^i})$ is exactly two points for
$i<\ell$, with cancelling local fixed point data, as we see by
applying the Lefschetz Fixed Point Formula and the $G$--Signature
Formula to
$T^{p^i}$.  Moreover, $\text{Fix}(T^{p^\ell}) - \text{Fix}(T) $ consists of a single orbit
of $p^\ell$ points, as one can see from the Lefschetz formula.

 Let $(a,b)$ and $(c,d)$
be the local fixed point data for $T$ (and hence $T^{p^\ell}$, too)
at the two points of $\text{Fix}(T)$.  Now the $G$--Signature Formula
for $T$ says that
$$
\frac{(\zeta^{a}+1)}{(\zeta^{a}-1)}\frac{(\zeta^{b}+1)}
{(\zeta^{b}-1)}+\frac{(\zeta^{c}+1)}{(\zeta^{c}-1)}\frac{(\zeta^{d}+1)}
{(\zeta^{d}-1)}=0
$$
where $\zeta=\exp(2\pi i/m)$.  By Theorem \ref{thm:t=0}, it follows
that the fixed point data cancels and we may assume that
$(c,d)=(a,-b)$, for example.  

Now let the fixed point data at
the $p^\ell$ new fixed points of
$T^{p^\ell}$ be
$(e,f)$.  Then the $G$--Signature Formula for $T^{p^\ell}$ says that
$$
p^\ell\frac{(\zeta^{e}+1)}{(\zeta^{e}-1)}\frac{(\zeta^{f}+1)}{(\zeta^{f}-1)}
=p^\ell$$
which cannot hold.  
\end{proof}

\begin{thm}[Restatement of Theorem \ref{thm:p=3}]\label{thm:p=3'}
There is no locally linear,\break
\emph{pseudofree} action of C$_9$ on
$\#_n\mathbf{C}P^2$ (or any closed, simply
connected, positive definite 4--manifold) with homology
representation of type
$$r(9)+s(3)+t(1)$$ where
$t\le 1$ and $s\ge 1$.
\end{thm} 
\begin{proof}
We will begin by analyzing a general
pseudofree action of $C_m$, where $m=p^2$, $p$ an odd prime.  We
will eventually get to a point where we can see that the cases
covered by Theorem \ref{thm:p=3} (\ref{thm:p=3'}) can be ruled out,
but that the door remains open for $p\ge 5$.

So suppose that the cyclic group $C_m$ has generator $T$ and acts on
$X=\#_n\mathbf{C}P^2$ with representation on $H_2$ given by
$$r(m)+s(p)+t(1)$$  There are two kinds of singular orbits:
\begin{enumerate}
\item $t+2$ fixed points $\{x_i\}$ with isotropy group
$C_m=\left<T\right>$.  The fact that $\text{card}X^{C_m}=t+2$
follows from the Lefschetz Fixed Point Formula.
\item $s$ orbits of size $p$, with isotropy group
$C_p=pC_m=\left<T^p\right>$:
$$\{y_i,Ty_i,\dots,T^{p-1}y_i:i=1,\dots,s\}$$ Again, the fact that
the number of such orbits is $s$  follows from the Lefschetz Fixed
Point Formula.
\end{enumerate}

We now interpret the $G$--Signature Formula in this situation.  On the
one hand we have $$\sigma(T,X)=t$$ and
$$\sigma(T^p,X)=t+ps.$$
Now let $T$ have local representation of type $(a_i,b_i)$ at
$x_i$.  Set
$\zeta=\text{exp}(2\pi i/m)$.  Then we have:
$$ t = \sigma(T,X) =
\sum_{i=1}^{t+2}\frac{(\zeta^{a_i}+1)}{(\zeta^{a_i}-1)}\frac{(\zeta^{b_i}+1)}
{(\zeta^{b_i}-1)}
$$ There are plenty of solutions to this equation coming from
equivariant connected sums of standard linear actions on
$\mathbf{C}P^2$. 

Similarly, let
$T^p$ have local type
$(c_i,d_i)$ (mod
$p$ rotation numbers) at
$y_i$.  Then the type of
$T^p$ at $T^jy_i$ is also $(c_i,d_i)$.  In addition, the type of
$x_i$ is again $(a_i,b_i)$, and the $G$--Signature Formula then yields
$$ t+ps = \sigma(T^p,X) =
\sum_{i=1}^{t+2}\frac{(\xi^{a_i}+1)}{(\xi^{a_i}-1)}\frac{(\xi^{b_i}+1)}
{(\xi^{b_i}-1)} + 
\sum_{i=1}^{s}p\frac{(\xi^{c_i}+1)}{(\xi^{c_i}-1)}\frac{(\xi^{d_i}+1)}
{(\xi^{d_i}-1)}
$$ where $\xi=\exp(2\pi i/p)=\zeta^p$.  One would like to think that
the two parts on each side of this equation correspond in the
``obvious'' way.

Under the substitution $\zeta\to\zeta^p$ we would obtain:
$$ t = \sigma(T,X) =
\sum_{i=1}^{t+2}\frac{(\xi^{a_i}+1)}{(\xi^{a_i}-1)}\frac{(\xi^{b_i}+1)}
{(\xi^{b_i}-1)}
$$ 
By Theorem \ref{thm:391} this holds for $p=3$ and $t\le 1$, but
does not hold in general. 
Whenever it does hold, however, we conclude that:
$$ s = 
\sum_{i=1}^{s}\frac{(\xi^{c_i}+1)}{(\xi^{c_i}-1)}\frac{(\xi^{d_i}+1)}
{(\xi^{d_i}-1)}
$$ 
For $p=3$ this last equation has no solutions, since the terms in
the sum are all $ \pm \frac{1}{3}$, and we are done.
\end{proof}

Even for $p>3$, if one can get to the last step of the proof, this 
formula does yield nontrivial restrictions.  One cannot then have
$s=1$, for example. 

\begin{rmk}At first glance one might think  that the equation
$$ s = 
\sum_{i=1}^{s}\frac{(\xi^{c_i}+1)}{(\xi^{c_i}-1)}\frac{(\xi^{d_i}+1)}
{(\xi^{d_i}-1)}
$$ 
 has no
solutions.  But in fact it does have solutions, at least for
appropriate $p$ and $s$:  Build an equivariant connected sum of
standard linear actions on $\#_{s} \mathbf{C}P^2$ (with $s+2$ fixed
points) until one sees a canceling pair of fixed point data. 
Removing that data yields the desired algebraic solution. This
procedure was used in earlier work of Edmonds and Ewing \cite{EdmondsEwing92} to
produce homologically trivial actions on a connected sum of copies
of $\mathbf{C}P^2$ with nonstandard fixed point data not coming from
an equivariant connected sum. Subsequently Hambleton and Lee \cite{Hambleton95} proved
that some of these locally linear actions cannot be smoothable.

The smallest solution and the only one we know for $s\le 2$ occurs
for $p=5$: 
$$
\frac{(\zeta^{1}+1)}{(\zeta^{1}-1)}\frac{(\zeta^{4}+1)}
{(\zeta^{4}-1)}+\frac{(\zeta^{2}+1)}{(\zeta^{2}-1)}\frac{(\zeta^{3}+1)}
{(\zeta^{3}-1)}=2
$$
where $\zeta=\exp(2\pi i/5)$, which we will use in the next section.
\end{rmk}

We now turn our attention to the actual construction
of an exotic, pseudofree, but not semifree, action of $C_{25}$, in
the case when $t=0$ and $s=2$.

\section{An exotic $C_{25}$ action} Here we give the
details of the construction of a pseudofree, but not
semifree, periodic map on a connected sum of copies of
$\mathbf{C}P^2$.

\begin{thm}[Restatement of Theorem \ref{thm:p=5}] There is
a locally linear {pseudofree} action of C$_{25}$ on
$\#_{10}\mathbf{C}P^2$ with homology representation of type
$2(5)$.
\end{thm} 

\begin{proof} Based on the work in the
preceding section we have an excellent picture of
what such an action would have to look like.  Let
$T$ denote a generator of
$C_{25}$.  Then
$T$ will have two fixed points, with canceling fixed point
data.   The fifth power $T^5$ will have two groups of five
fixed points, in addition to the two original fixed points
of $T$.  Each of the five points in each group must have
the same fixed point type, since the points are permuted
cyclically by $T$.

   We will use $(1,1), (1,24)$ for the fixed point
data associated to the two fixed points of $T$.  And
we will have one orbit of five points all of which
have $(1,4)$ as their local representation of $T^5$
(mod 5); and a second orbit of five points all of
which have ($2,3)$ as their local representation of
$T^5$ (mod 5).  Thus, for the desired action on the
expected 4--manifold we will have
$\sigma(T,M^4)=0$ and $\sigma(T^5,M^4)=10$.

We will start with a 4-ball $D^4$ on which $T$ acts with
local representation of type $(1,1)$.  We will then add two
orbits of five 2--handles each along a carefully chosen
equivariantly framed link.  Each 2--handle will be invariant
under T$^5$ and will contribute one isolated fixed point of
$T^5$.  The integral framings are determined mod 5 by the
required local fixed point data.  Our challenge is to find
the appropriate 10 by 10 equivariant integral linking
matrix with the additional requirement that the
corresponding integral bilinear form is equivalent to the
standard $10\left<+1\right>$.

For reference we record three lemmas from Edmonds and Ewing \cite{EdmondsEwing92}.  For
the lemmas suppose the generator $T$ of a free action of $C_m$ acts
on
$D^4=D^4(a,b)$ by
$T(z,w)=(\zeta^az,\zeta^bw)$, where $\zeta=\exp(2\pi
i/m)$ and $a$ and
$b$ are prime to $m$.  

\begin{lem}
If $S_x$ and $S_y$ are two disjoint, invariant, simple
closed curves in $\partial D^4(a,b)$ on which the
generator $T$ of $C_m$ operates by rotation by $2\pi x/m$
and $2\pi y/m$, respectively, then $Lk(S_x,S_y)\equiv
abx'y'\mod m$, where $xx'\equiv 1\mod m$ and $yy'\equiv
1\mod m$.
\end{lem}

Recall the standard observation that in the orbit
space $\partial D^4/C_m$ the simple closed curve
$\overline{S}_x$ represents 
$T^{x'}$ of $\pi_1(\partial D^4/C_m)$ under the
standard identification with $C_m=\left<T\right>$,
where $xx'\equiv 1\mod m$.

\begin{lem}\label{lem:stddata}
Let $f\co S^1\times D^2\to S^1\times D^2$ by
$f(z,w)=(z,z^rw)$.  If $C_m$ acts on the target
$S^1\times D^2$ by $T(z,w)=(\zeta^az,\zeta^bw)$, then $f$
is equivariant with respect to the action on the domain
given by $T(z,w)=(\zeta^az,\zeta^{b-ra}w)$.
\end{lem}

\begin{lem}
Suppose a 2--handle is added to $D^4(a,b)$ along a
curve $S_k\subset \partial D^4(a,b)$, with framing
$r$.  Then the framing map can be chosen so that
the action of $C_m$ extends over the 2--handle
$B^2\times D^2$ with action of type
$(-k,-rk+abk') \mod m$.
\end{lem}

Now back to the situation of $m=25$ and of $C_{25}$
acting on D$^4(1,1)$. To add a 2--handle invariant
under
$T^5$ that gives rise to a fixed point of $T^5$ of
type $(1,4)$ we need a curve $S_{-1}=S_4$ on which
$T^5$ acts by rotation by
$-2\pi/5$ and with framing $r$ so that
$-r*(-1)+1*1*(-1)\equiv 4\mod 5$, and hence $r\equiv
0\mod 5$.  We will use five such curves, cyclically
permuted by $T$.

Similarly, to add a 2--handle invariant under $T^5$ that
gives rise to a fixed point of $T^5$ of type $(2,3)$ we
need a curve $S_{-2}=S_3$ on which $T^5$ acts by rotation
by $-2*2\pi/5$ and with framing $r$ so that
$-r*(-2)+1*1*2\equiv 3\mod 5$, and hence $r\equiv
3\mod 5$.  We will also use five such curves, again
cyclically permuted by $T$.

We obtain these curves by representing $T^{20}$ and
$T^{10}$, respectively, by disjoint simple closed
curves in
$L=\partial D^4(1,1)/C_{25}$ using the canonical
identification of $\pi_1(L)$ with $C_{25}$.

Attaching two $2$--handles downstairs, or ten 2--handles
upstairs, we end up with a 2 by 2 linking matrix over
$\mathbf{Z}[C_{25}]$.  If $W$ is the resulting 4--manifold,
then
$H_2(W;\mathbf{Z}[C_{25}])=H_2(\widetilde{W};\mathbf{Z})=
\mathbf{Z}[C_5]\oplus \mathbf{Z}[C_5]$.  The
$\mathbf{Z}[C_{25}]$--valued linking form depends on the
exact choice of curves and integral framings and integral
linking numbers.  But it is determined mod 5 by the
above lemmas as:
$$
\left[
\begin{matrix}
5+T^{5}+T^{10}+T^{15}+T^{20} &
3+3T^{5}+3T^{10}+3T^{15}+3T^{20} \\
3+3T^{5}+3T^{10}+3T^{15}+3T^{20} &
3+4T^{5}+4T^{10}+4T^{15}+4T^{20}
\end{matrix}
\right]
$$
Moreover, having made one tentative choice of these
curves and their framings, one can make some changes
by sliding a small arc on one curve across another of
the curves and extending equivariantly.  The upshot
of this is that one can achieve any integral matrix
of the form
$$
\left[
\begin{matrix}
a_0+a_1(T^{5}+T^{20})+a_2(T^{15}+T^{10}) &
b_0+b_1(T^{5}+T^{20})+b_2(T^{15}+T^{10}) \\
b_0+b_1(T^{5}+T^{20})+b_2(T^{15}+T^{10})  &
c_0+c_1(T^{5}+T^{20})+c_2(T^{15}+T^{10})
\end{matrix}
\right]
$$
where $a_0\equiv 0\mod 5$, $a_1\equiv a_2\equiv 4\mod
5$, $b_0\equiv b_1\equiv b_2\equiv 3\mod 5$, 
$c_0\equiv 3\mod 5$,
$c_1\equiv c_2\equiv 4\mod 5$.
The challenge is to choose the integers
$a_0,\ldots,c_2$ so that the resulting underlying
rank 10 integral form is unimodular and in fact
represents $10\left<+1\right>$.

It took some time to find such a matrix, by a combination
of brute force computer searching and judicious
simplifications.  Here is one such matrix:
$$
\left[
\begin{matrix}
45+16(T^{5}+T^{20})+16(T^{15}+T^{10}) &
93+23(T^{5}+T^{20})+23(T^{15}+T^{10}) \\
93+23(T^{5}+T^{20})+23(T^{15}+T^{10}) &
198+29(T^{5}+T^{20})+29(T^{15}+T^{10})
\end{matrix}
\right]
$$
One can check that this matrix has determinant 1 in
$\mathbf{Z}[C_{25}]$\ !

 The corresponding
integer-valued form is given by:
$$
\left[
\begin{matrix}
45&16&16&16&16&93&23&23&23&23\\
16&45&16&16&16&23&93&23&23&23\\
16&16&45&16&16&23&23&93&23&23\\
16&16&16&45&16&23&23&23&93&23\\
16&16&16&16&45&23&23&23&23&93\\
93&23&23&23&23&198&29&29&29&29\\
23&93&23&23&23&29&198&29&29&29\\
23&23&93&23&23&29&29&198&29&29\\
23&23&23&93&23&29&29&29&198&29\\
23&23&23&23&93&29&29&29&29&198\\
\end{matrix}
\right]
$$
One can check directly that this matrix has
determinant 1 in $\mathbf{Z}$, is positive definite,
and in fact can be reduced to a $10\times 10$
identity matrix by simultaneous integral row and
column operations.  Therefore it represents the
standard intersection form $10\left<+1\right>$ (as
opposed to $E_8\oplus 2\left<+1\right>$).  It helps
to use a computer algebra package, not only to find
the matrix, but even just to verify its claimed
properties. Alternatively one might be able to take a more theoretical approach, showing that there exist three independent vectors of norm $1$ and invoking the known classification of unimodular forms of low rank. Or one might be able to apply N. Elkies's \cite{Elkies95} characterization of the standard lattice by showing that the present lattice has no ``short characteristic vectors''. We have, however, opted for a more direct and na\"ive approach. See the end of this section for the Maple code used to make these direct verifications.

Form a smooth 4--manifold $W^4$ with smooth $C_{25}$
action by attaching these ten 2--handles to $D^4(1,1)$
equivariantly.  The matrix above gives the
intersection pairing on $H_2(W)$.  Because the
intersection pairing is unimodular, the boundary 
$\Sigma$ of $W^4$ is an integral homology 3-sphere. 
According to Freedman's fundamental work \cite{Freedman82}, such a
homology 3-sphere is the boundary of a compact
contractible, topological 4--manifold $\Delta$.  We
must extend the action on $\Sigma$ to a locally
linear action on $\Delta$.  This problem was
completely analyzed in Edmonds \cite{Edmonds87} for $m$ prime, and
then in Kwasik and Lawson \cite{KwasikLawson93} for general $m$.  The
necessary and sufficient conditions are that the
quotient $Q=\Sigma/C_m$ have the same signature
invariants ($\rho$ or $\alpha$ invariants) as some
lens space, and the  same Reidemeister torsion, up
to squares, of the  same lens space.

Our choice of fixed point data guarantees that the
signature invariants are correct.

The torsion condition comes down to the  issue of
the determinant of the above matrix over
$\mathbf{Z}[C_{25}]$ being a square of a unit in
$\mathbf{Z}[C_{25}]/(N)$, where $N\in
\mathbf{Z}[C_{25}]$ denotes the sum of the group
elements.  See the argument in Edmonds and Ewing \cite{EdmondsEwing92}, 
pages 1115-1117.  We have the good fortune that this
determinant is in fact 1 in
$\mathbf{Z}[C_{25}]$.
\end{proof}

\begin{rmk}
What about the prospects of similar actions for
$C_m$, $m=p^2$, $p> 5$?  We can find suitable
fixed point data, although the number of terms
necessary seems to grow with $p$.  The problem of
lifting the data to find a suitable positive definite
linking matrix, however, seems formidable at the
present time.  In addition it is unclear whether the
torsion issue would be so easily resolvable in
general.
\end{rmk}
\subsection*{Maple Code for Matrix Reduction}\label{maple}
The following Maple code shows that the $10\times 10$ matrix discussed above can be reduced over the integers to the identity matrix by a suitable sequence of 78 simultaneous row and column operations.
{\small
\begin{verbatim}
> with(linalg):
> m:=matrix([
> [45,16,16,16,16,93,23,23,23,23],
> [16,45,16,16,16,23,93,23,23,23],
> [16,16,45,16,16,23,23,93,23,23],
> [16,16,16,45,16,23,23,23,93,23],
> [16,16,16,16,45,23,23,23,23,93],
> [93,23,23,23,23,198,29,29,29,29],
> [23,93,23,23,23,29,198,29,29,29],
> [23,23,93,23,23,29,29,198,29,29],
> [23,23,23,93,23,29,29,29,198,29],
> [23,23,23,23,93,29,29,29,29,198]]):
> # Procedure to test whether a "local move" can reduce the 
> # (sum of the) diagonal entries
> test:=proc(a) local i,j,n,tt;
> tt:=0;
> n:=rowdim(a);
> for i from 1 to n-1 while tt=0 do
> for j from i+1 to n while tt=0 do
> if 2*abs(a[i,j])>min(abs(a[i,i]),abs(a[j,j])) then tt:=1 fi;
> od;
> od;
> tt
> end:
> # Procedure to perform local moves as long as they can reduce the 
> # (sum of the) diagonal entries
> reduce:=proc(a) local i,j,n,b,e,k;
> global count;
> b:=a;
> n:=rowdim(b);
> if test(evalm(b))=0 then RETURN(evalm(b)) fi;
> for i from 1 to n do
> if b[i,i]<>0 then
> for j from 1 to n do
> if j<>i then
> if abs(2*b[i,j])>abs(b[i,i])  
> then 
> count:=count+1;
> k:=floor(abs(b[i,j])/abs(b[i,i]));
> e:=signum(b[i,j]/b[i,i]);
> b:=addrow(addcol(b,i,j,-e*k),i,j,-e*k); 
> fi; fi; od; fi; od;
> if test(evalm(b))=0 then RETURN(evalm(b)) else
> RETURN(evalm(reduce(b))) fi
> end:
> #Run the routines on the given matrix
> count:=0:
> reduce(m);
> print(number_of_simultaneous_row_and_column_operations=count):
\end{verbatim}
}

\end{document}